\documentclass{article}
\usepackage{amsmath,amssymb,indentfirst}
\usepackage[latin1]{inputenc}
\usepackage{epsfig}
\usepackage{psfrag}
\usepackage{hyperref}
\usepackage{url}
\usepackage{tikz}

\newtheorem{theorem}{Theorem}
\newtheorem{lemma}{Lemma}

\newtheorem{cor}{Corollary}

\newtheorem{exe}{Example}
\newtheorem{rmk}{Remark}
\usepackage{fancyhdr}
\usepackage{ifthen}

\newcommand{\qed}{\hfill\mbox{\raggedright $\Box$}\medskip}

\begin{document}

\title{THE CHANGE OF VARIABLE FORMULA INTEGRALS,\\ 
DO THEY HAVE EQUAL VALUE?}

\author{Oswaldo Rio Branco de Oliveira}
\date{}
\maketitle

\begin{abstract} Assuming that the two integrals in the {\sl{Change of Variable Formula
for the Riemann integral on the real line}} are finite,  one can rightfully
ask if they have equal value. We give a positive answer to this question. 
The proof is very easy to follow and to keep in mind. The result is extended to
improper Riemann integrals. Three examples are given.
\end{abstract}


\hspace{- 0,6 cm}{\sl Mathematics Subject Classification: 26A06, 26A42}

\hspace{- 0,6 cm}{\sl Key words and phrases:} One Variable Calculus, Integrals of Riemann.

\section{Introduction}


This article deals with the Riemann integrals

$$ \int_{\varphi(\alpha)}^{\varphi(\beta)}f(x)dx \  \textrm{and}\ \int_\alpha^\beta f(\varphi(t))\varphi'(t)dt, \leqno{(1)}$$ 
appearing on many versions of the {\sl{Change of Variable Theorem on the real line}}. If these two integrals are finite, it is natural to ask if their values are equal. In other words, we investigate if the existence of these integrals is already enough to guarantee that the {\sl{Change of Variable Theorem}} is valid. We give a positive answer to this question. Three examples are shown.

The very influential general version due to  Kestelman \cite{Kestelman}, see also Davies \cite{Davies}, establishes the integral identity
$$ \int_{G(\alpha)}^{G(\beta)}f(x)dx = \int_\alpha^\beta f(G(t))g(t)dt,$$
where $g: [\alpha,\beta] \to \mathbb R$ is integrable, the substitution map $G$ is given by
$$G(t)= G(\alpha) + \int_\alpha^t g(\tau)d\tau,\ \textrm{where} \ t \in [\alpha,\beta],$$
and the base function $f:G([\alpha,\beta]) \to \mathbb R$ is integrable.  It is worth emphasizing that Kestelman's version requires the integrability of $g$, the derivative of the map $G$. The identity $G'(t) = g(t)$ is true on all points of continuity of $g$.

Many articles have been written about Kestelman's version. See Bagby \cite{Bagby}, Kuleshov \cite{Kuleshov}, Pouso \cite{Pouso}, Sarkhel and V\' yborn\' y \cite{Sarkhel}, Tandra \cite{Tandra}, and Torchinsky \cite{Torchinsky1}. See also Torchinsky's book \cite{Torchinsky2}.

Regarding the change of variable theorem with the integrals having the format described in (1), de Oliveira \cite{Oliveira} shows a very simple change of variable formula that we improve a bit on section 4 and then use in section 5 of the present article. Such change of variable, suitable for first year Calculus, has no condition on the derivative of the substitution map $\varphi$ other than the finiteness of the integral of the product $(f\circ \varphi)\varphi'$ but on the other hand requires the existence of a primitive of the function $f$. Nevertheless,  this elementary change of variable is quite useful to broaden the scope of the {\sl Change of Variable Formula}.

We present three other versions for the {\sl{Change of Variable Theorem for the Riemann integral on the real line}}. These three have the following features:
\begin{itemize} 
\item[$\bullet$] Theorem \ref{T1} already covers practical situations, when one deals with a concrete and explicit function $f$. Its proof  uses a minimum of {\sl{ Lebesgue's Integration Theory}}, but not an actual Lebesgue integral.  We comment more on these facets on the final sections  {\sf{Remarks}} and  {\sf{ Conclusion}}.
\item[$\bullet$] The proof of Theorem \ref{T2} is an adaptation of that of Theorem \ref{T1}. However, it relies a little bit more on {\sl{Lebesgue's Integration}}. This version is general and puts no condition on the derivative of the substitution map $\varphi$ other than the finiteness of the integral of $(f\circ \varphi)\varphi'$. Moreover, following a remark made by Kuleshov \cite{Kuleshov}, we argue that the base function $f$ can be unbounded outside its interval of integration (see Remark \ref{R4}). Theorem \ref{T2} implies a Change of Variable Formula for improper Riemann integrals.
\item[$\bullet$] Theorem \ref{T3} is immediate from Theorem \ref{T2} and has a very trivial statement.
\end{itemize}

\section{Notation}

Let $[a,b]$ be a closed and bounded interval on the real line. Given an arbitrary bounded function $f:[a,b] \to \mathbb R$, its  {\sl lower Darboux sum} is given by
 $$s(f,\mathcal{P}) = \sum _{i=1}^n m_i\Delta x_i,$$
 where $\mathcal{P} = \{a=x_0\leq x_1\leq \cdots \leq x_n =  b\}$ is a {\sl{ partition}} of the interval $[a,b]$, $m_i = \inf \{f(x): x \in [x_{i-1}, x_i]\}$, and $\Delta x_i = x_i - x_{i-1}$, for each $i=1,\ldots, n$. Moreover, the {\sl{norm}} of the partition $\mathcal{P}$ is given by $|\mathcal{P}| = \max \{\Delta x_1, \ldots, \Delta x_n\}$. 

Given a Riemann integrable function $f:[a, b] \to \mathbb R$, it is well known that
$$\lim_{|\mathcal{P}|\to 0} s(f,\mathcal{P}) = \int_a^b f(x)\,dx.$$

We denote the Lebesgue measure on the real line by $m$. Given a Lebesgue integrable function $f:[a,b]\to \mathbb R$, we denote its Lebesgue integral by
$$\int_a^bf(x)dm.$$ 
It is well known that if a function $f:[a,b] \to \mathbb R$ is Riemann integrable, then $f$ is Lebesgue integrable and both integrals have equal value.

Unless specified otherwise, all the integrals that follow are Riemann integrals.

Given an arbitrary function $f:I \to \mathbb R$, with $I$ an interval on the real line, we say that $F:I \to \mathbb R$ is a {\sl primitive} of $f$ if $F'=f$. At the possible endpoints, the derivative $F'$ is understood as a lateral derivative.

We say that a subset $\mathcal{N}\subset \mathbb{R}$ is a {\sl{ null}} set, or that it has {\sl{ measure zero}}, if for every $\epsilon >0$ there exists a sequence of open intervals $I_n=(a_n,b_n)$ on the real line such that $\mathcal{N}$ is contained on the union of these open intervals $I_{n's}$, with $[(b_1-a_1) + (b_2 - a_2) + \cdots]<\epsilon$.

It is not hard to see that a countable union of null sets is also a null set.

Let us now consider a finite collection $\{[x_i, x_{i+1}]: 1 \leq i \leq n\}$ of consecutive bounded and closed intervals on the real line . A real function of one real variable $g: [x_1, x_2]\cup [x_2,x_3]\cup \cdots \cup [x_{n},x_{n +1}] \to \mathbb R$ is   {\sl{ piecewise linear continuous}} if $g$ is continuous and the restriction of $g$ to each sub-interval $[x_i,x_{i+1}]$ is an affine function. With some abuse of terminology, we also call this restriction a {\sl linear function}. Clearly, two consecutive sub-intervals  have one common point.

Given two  points in the plane,  $(x_1,y_1)$ and $(x_2, y_2)$, with $x_1 < x_2$, there exists {\sl{a unique linear continuous function}}  $g:[x_1, x_2]\to \mathbb R$ whose graph connects  the points $(x_1,y_1)$ and $(x_x,y_2)$. Evidently, the function $g$ is a constant if $y_1=y_2$, linear increasing if $y_1<y_2$, and linear decreasing if $y_1>y_2$.

 %

\section{An approximation lemma}

In this section we do not employ any Lebesgue's measure theory. 

In what follows, fixed an arbitrary non-negative function $f: [a,b] \to [0,+\infty)$ 
we define an infinite sequence of non-negative piecewise linear continuous functions $f_n: [a,b] \to [0,+\infty)$ that approximates the function $f$ from below. 

{\sl{ The definition of $f_1$ and $f_2$}}. We define the functions $f_1: [a,b] \to [0,+\infty)$ and $f_2:[a,b]\to [0, + \infty)$ as $f_1=f_2=0$, the zero function.

{\sl{ The definition of $f_n$, for $n\geq 3$}}. We start by dividing $[a,b]$ into $2^{n}$ non-overlapping closed sub-intervals of equal length $(b -a)/2^n$. Let us label these sub-intervals as $I_1, I_2, \ldots, I_{2^{n}}$, naturally ordered from left to right along the real line. 
The final endpoint of $I_k$ is the initial endpoint of $I_{k+1}$, for all possible $k$.

Next, we divide each sub-interval $I_k$ into four  non-overlapping closed sub-intervals labeled as $I_k^1$, $I_k^2$, $I_k^3$, and $I_k^4$ also naturally ordered from left to right along the real line. For these four, we choose their lengths conveniently.

{\sl The case $ k \in \{2,\ldots, 2^{n} -1\}$}. The sub-intervals $I_k^1$ and $I_k^4$ have equal length 
$$\epsilon_n = \frac{b -a}{n2^{n}}. \leqno (2)$$ 
The sub-intervals $I_k^2$ and $I_k^3$ (these are ``on the middle'' of $I_k$) have equal length. Evidently, the sum of the lengths of these four sub-intervals is $(b-a)/2^n$.

{\sl The case $k=1$}. Then,  the sub-interval $I_1^4$  has length $\epsilon_n$. The other three, $I_1^1$, $I_1^2$, and $I_1^3$, have equal length. 
The total sum of their lengths is $(b-a)/2^n$.

{\sl The case $k=2^{n}$}. Then,  the sub-interval $I_{2^{n}}^1$ has length $\epsilon_n$. The other three, $I_{2^{n}}^2$, $I_{2^{n}}^3$, and $I_{2^{n}}^4$, have equal length. 
The total sum of their lengths is $(b-a)/2^n$.

Next, let $m_k$ be the infimum of $f$ restricted to $I_k$, for each $k=1,\ldots, 2^{n}$. 

We define $f_n$ as a constant in the following sub-intervals,
$$f_n = 
\left\{\begin{array}{lll}
m_1,\ \textrm{on the interval} \ I_1^1\cup I_1^2 \cup I_1^3\\
m_k, \ \textrm{on the interval}\ I_k^2\cup I_k^3,\ \textrm{for each} \ k \in \{2,\ldots, 2^{n} -1\},\\
m_{2^n}, \ \textrm{on the interval} \ I_{2^{n}}^2 \cup I_{2^{n}}^3\cup I_{2^{n}}^4.
\end{array}
\right.
$$

Let us define $f_n$ on the ``level change'' intervals $I_k^4$, $I_{k+1}^1$, for $k = 1,\ldots, 2^{n} -1$. 

If $m_k = m_{k+1}$, we just put $f_n =m_k$ on $I_k^4\cup I_{k+1}^1$.

If $m_k < m_{k+1}$, we put $f_n = m_k$ on $I_k^4$. At this juncture, $f_n$ is defined at both endpoints of the closed $I_{k+1}^1$ (this intersects $I_k^4$ and $I_{k+1}^2$) and thus we may define $f_n$ as linear, continuous, and increasing on  $I_{k+1}^1$ (see Figure 1). We have the inequalities $m_k \leq f_n(x)\leq m_{k+1}$ for all  $x$ in the ``level change'' interval $I_{k+1}^1$.
\begin{figure}[h]
  \centering
  \vspace{-0,4 cm}
  \begin{tikzpicture}[scale=0.76]
    \draw[-stealth] (-0.5,0)--(6,0) node[below]{$x$};
    \draw[-stealth] (0,-0.5)--(0,4) node[right]{$y$};
    \draw[dashed] (0,1.5)node[left]{$m_k$} -- (0.1,1.5);
    \draw[dashed] (0,3)node[left]{$m_{k+1}$} -- (0.1,3) (3,3) -- (3.5,3);
		\draw[thick] (1,1.5)--(3,1.5)--(3.5,3)--(5,3);
		\draw[dashed] (1,0)--(1,1.5);
    \draw[dashed] (3,0.1)--(3,0);
    \draw[dashed] (3.5,0.1)--(3.5,0);
    \draw[dashed] (5,0)--(5,3);
		\draw[dashed] (3,1.5)--(3,3);
    \node[below] (Ik) at (3.3,0) {$I_{k+1}^1$};
  \end{tikzpicture}
  \label{fig:fig1}
  \vspace{-0,4 cm}
  \caption[Short]{Graph of $f_n$ over $(I_k^2 \cup I_k^3 \cup I_k^4) \cup I_{k+1}^1 \cup (I_{k+1}^2\cup I_{k+1}^3)$, with
    $m_k<m_{k+1}$.}
		\vspace{-0,25 cm}
\end{figure}

If $m_k > m_{k+1}$, we put $f_n = m_{k+1}$ on $I_{k+1}^1$. Similarly, by now $f_n$ is defined at both endpoints of the closed interval $I_{k}^4$ (which intersects $I_k^3$ and $I_{k+1}^1$) and we define $f_n$ as linear, continuous, and decreasing on $I_{k}^4$ (see Figure 2). We have the inequalities $m_{k+1}\leq f_n(x) \leq m_k$ for all $x$ in the ``level change'' interval $I_k^4$.
\begin{figure}[h]
  \centering
	\vspace{-0,2 cm}
  \begin{tikzpicture}[scale=0.76]
    \draw[-stealth] (-0.5,0)--(6,0) node[below]{$x$};
    \draw[-stealth] (0,-0.5)--(0,4) node[right]{$y$};
    \draw[dashed] (0,3)node[left]{$m_k$} -- (0.1,3) (2.5,3) -- (3,3);
    \draw[dashed] (0,1.5)node[left]{$m_{k+1}$} -- (0.1,1.5);
		\draw[thick] (1,3)--(2.5,3)--(3,1.5)--(5,1.5);
		\draw[dashed] (1,0)--(1,3);
    \draw[dashed] (2.5,0.1)--(2.5,0.0);
		\draw[dashed] (3,0.1)--(3,0.0);
    \draw[dashed] (3,2)--(3,3);
    \draw[dashed] (5,0)--(5,1.5);
    \node[below] (Ik) at (2.7,0) {$I_{k}^4$};
  \end{tikzpicture}
  \label{fig:fig2}
  \vspace{-0,4 cm}
  \caption[Short]{Graph of $f_n$ over $(I_k^2 \cup I_k^3) \cup I_k^4 \cup (I_{k+1}^1 \cup I_{k+1}^2\cup I_{k+1}^3)$, with $m_k>m_{k+1}$.}
\end{figure}

\vspace{- 0,2 cm}
All in all, we have  $\min(m_k, m_{k+1}) \leq f_n(x) \leq \max(m_k, m_{k+1})$ if $x\in I_k^4 \cup I_{k+1}^1$.

The definition of the sequence is finished. 
Two remarks are appropriate.

First, the graph of the restriction $f_n: I_k \to [0, +\infty)$ is below the graph of the constant function $ m_k: I_k \to [0,+\infty)$,  for each $k$. 

Second, $f_n:[a,b] \to [0,+\infty)$ is continuous and  satisfies $0\leq f_n \leq f$. 

\begin{lemma} \label{L1}{\sl (Approximation Lemma.)} The sequence of piecewise linear functions $f_n:[a,b]\to [0,+\infty) $ defined right above has the following four properties.
\begin{itemize}
\item[1.] We have $0\leq f_n \leq f$, for all $n$.
\item[2.] Each $f_n$ is a continuous function.
\item[3.] If $f$ is continuous at a point $p$,  then $ \lim f_n(p)=f(p)$.
\item[4.] If $f$ is integrable on the sub-interval $[c,d]$, then we have
$$\int_c^d |f_n(x) - f(x)|dx \xrightarrow[\raisebox{0.8ex}{\scriptsize $n\!\to\! +\infty$}]{} 0.$$
\end{itemize}
\end{lemma}
{\bf Proof.} The properties 1 and 2 are immediate. Let us show 3 and 4.
\begin{itemize}
\item[3.] {\sl{ The case $p \in (a,b)$.}} Given an arbitrary $\epsilon >0$, there exists a small enough $\delta >0$ such that  $f(x) \in [f(p) - \epsilon, f(p) + \epsilon]$ for all $x \in [p- \delta, p + \delta]\subset (a,b)$. Let us pick $n$ big enough (we may assume $n>3$) such that 
$$\frac{b-a}{n} < \delta/4$$

Then, following the definition of $f_n$, its domain $[a,b]$ is subdivided into $2^{n}$ closed sub-intervals of length $(b-a)/2^{n}$. Among these sub-intervals, there exists a  $I_k$ such that 
$p \in I_k$.  We have $2^{n} >n$, which implies that the length of the interval $I_k$ is smaller than $\delta/4$. Thus,  $p \in I_k$ and we definitely have the inclusions
$I_{k-1}\cup I_k\cup I_{k+1} \subset [p - \delta, p + \delta] \subset (a,b)$. Hence, $f(x) \in [f(p)-\epsilon, f(p) + \epsilon]$ for every $x\in I_{k-1}\cup I_k\cup I_{k+1}$. Following the notation on the definition of the sequence $(f_n)$, we conclude that the infimums $m_{k-1}, m_k$, and $m_{k+1}$ are all in the interval $[f(p) - \epsilon, f(p) + \epsilon]$. Also by the construction of the sequence $(f_n)$, we have the inequalities
$$\min \{m_{k-1},m_k,m_{k+1}\} \leq f_n(x)\leq \max \{m_{k-1},m_k, m_{k+1}\}\, \textrm{for all} \ x \in I_k.$$
In particular, we find that $f_n(p) \in [f(p) - \epsilon, f(p) + \epsilon]$, for all $n$ big enough. This shows that  $f_n(p) \to f(p)$ as $n \to +\infty$.

{\sl{ The case $p=a$}}, assuming $f$ continuous at $p=a$. Given $\epsilon >0$, there exists $\delta >0$ such that $f([a,a + \delta]) \subset [f(a) - \epsilon, f(a) + \epsilon]$. Then, given $n$ big enough,  $[a,b]$ is divided into sub-intervals $I_1, \ldots, I_{2^{n}}$ of equal length, with $a \in I_1  \subset [a,a+\delta]$. By  definition  we have
$f_n(a) = m_1$, with $m_1$ the infimum of $f$ on $I_1^1 \cup I_1^2\cup I_1^3 \subset I_1 \subset [a,a+\delta]$. Hence, $f_n(a) = m_1 \in [f(a) - \epsilon, f(a) + \epsilon]$.

{\sl{ The case $p=b$}} is analogous to the case $p=a$.

\item[4.] Let $M$ be a constant satisfying  $0\leq f(x)\leq M$, for all points $x\in [c,d]$. Given $n$, we consider the sub-intervals $I_k$ of  $[a,b]$,  where $1\leq k\leq 2^{n}$, as  in the construction of  
$f_n$. Next, let $\mathcal{P}_n$ be the partition of the interval $[c,d]$ determined by all the non-empty intersections $I_k\cap[c,d]$, where $1\leq k\leq 2^{n}$.
As previously remarked, the graph of the restriction $f_n: I_k \to [0,+\infty)$ is below the graph of the constant function $m_k: I_k \to [0,+\infty)$, for each $k$. 
 Thus, with some abuse of notation, the lower Darboux sum $s(f,\mathcal{P}_n)$ of the restriction $f:[c,d]\to \mathbb R$ and the Riemann integral of the continuous restriction $f_n: [c,d]\to \mathbb R$ satisfy the following inequalities [right at this juncture we take advantage of the choice $\epsilon_n = (b-a)/n2^{n}$, see equation (2), we also use that the restriction $f:[c,d]\to \mathbb R$ is bounded by $M$,  define the value of an integral over the empty set as zero, and remark that there are at most $2^{n-1}$ ``level change'' intervals (see Figure 1 and Figure 2), all with length $\epsilon_n$, on the construction of the piecewise linear function $f_n$],
\begin{eqnarray*}
(3)\ \ \ \ \,   0 \leq s(f,\mathcal{P}_n) - \int_c^d f_n(x)dx & =  &\sum_{k=1}^{2^n} \int_{I_k\cap [c,d]} [m_k - f_n(x)]dx\\
                                              &\leq& \sum_{k=2}^{2^{n}} \frac{\epsilon_n|m_k - m_{k-1}|}{2}\leq \frac{M(b-a)}{n}. \  \ \ \ \, \\
\end{eqnarray*}
Obviously, the norm of the partition $\mathcal{P}_n$ goes to $0$ as $n$ goes to $+\infty$.

Thus, since the restriction $f:[c,d] \to \mathbb R$ is integrable, its lower Darboux sums $s(f,\mathcal{P}_n)$ converge to its integral as $n$ goes to $+\infty$. 
From this, and from the inequalities (3), it follows that 
$$ \int_c^d f_n(x)dx \xrightarrow[\raisebox{0.8ex}{\scriptsize $n\!\to\! +\infty$}]{}  \int_c^d f(x)dx.$$
 Therefore, from the inequalities  $0\leq f_n \leq f$ it follows  that the integrals of $|f_n-f|=f-f_n$,  over the interval $[c,d]$, converge to $0$ as $n$ goes to $+\infty$.
\end{itemize}
The proof of the lemma is complete.$\qed$

\section{ Theorems used}

Besides assuming the well-known {\sl{ Lebesgue's Dominated Convergence Theorem}}, and enunciating an elementary version of {\sl{ The Fundamental Theorem of Calculus for the Lebesgue Integral}},  we enunciate and use four other lemmas. We prove two of these lemmas. One of these two, is a very basic lemma in Differential Calculus. The last two results in this section, Lemma \ref{openness}  and Lemma \ref{Saks}, are not used in the proofs of any the three theorems that form the center of this article (namely, Theorem \ref{T1}, Theorem \ref{T2}, and Theorem \ref {T3}). We only employ them to prove a corollary of Theorem \ref{T2} (namely, Corollary \ref{C2}).

We begin by the following quite intuitive result. 

\begin{lemma} \label{Serrin} ({\sf Serrin and Varberg's Theorem.}) Let $\varphi:\mathbb R\to \mathbb R$ be an arbitrary function.  Let us suppose that $\varphi$ has a derivative (finite or infinite) on a set $\mathcal{M}$ such that $\varphi(\mathcal{M})$ is a null set. Then, we have $ \varphi' = 0$ almost everywhere on $\mathcal{M}$. 
\end{lemma}
{\bf Proof.} See Serrin and Varberg \cite{Serrin}. Their proof is nice, short, and easy.$\qed$


%

%

\begin{lemma} \label{Oliveira} ({\sf An elementary Change of Variable Theorem for the Riemann integral.}) Let us consider $f:I \longrightarrow \mathbb R$, where $I$ is an interval and $f$ has a primitive, and a continuous map $\varphi:[\alpha,\beta]\longrightarrow I$ that is differentiable on $(\alpha,\beta)$. Let $J$ be the closed interval with endpoints $\varphi(\alpha)$ and $\varphi(\beta)$. The following is true.  
\begin{itemize}
\item[$\bullet$] If the product $(f\circ \varphi)\varphi'$ is integrable and $f$ is integrable on $J$, then we have 
$$\int_{\varphi(\alpha)}^{\varphi(\beta)} f(x)dx = \int_\alpha^\beta f(\varphi(t))\varphi'(t)dt.$$
\end{itemize}
\end{lemma}
{\bf  Proof.} (This lemma  takes into consideration a remark made by Kuleshov \cite{Kuleshov},  
and improves a version of the change of variable theorem  by de Oliveira \cite{Oliveira}.)

If the one-sided derivatives of the substitution map $\varphi$ at the endpoints $\alpha$ and $\beta$ do not exist, then we may define them at will. This does not affect the existence or the value of the integral of the product $(f\circ\varphi)\varphi'$.

Now, let $F:I \to \mathbb R$ be a primitive of $f$ and consider a small enough $\epsilon>0$. We have the identities
\begin{eqnarray*}
\int_{\alpha + \epsilon}^{\beta - \epsilon} f(\varphi(t))\varphi'(t)dt &=& \int_{\alpha + \epsilon}^{\beta - \epsilon}(F\circ\varphi)'(t)dt\\
                                                                       &=& F(\varphi(\beta- \epsilon)) -F(\varphi(\alpha + \epsilon)).
\end{eqnarray*}
Letting $\epsilon \to 0$, by the continuity of the composite $F\circ \varphi$ and the continuity of the Riemann integral with respect to the integration endpoints we arrive at
$$ \int_{\alpha}^{\beta} f(\varphi(t))\varphi'(t)dt = F(\varphi(\beta)) -F(\varphi(\alpha)). $$
Since $F$ is a primitive of $f$, and the function $f$ is Riemann integrable on the interval with endpoints $\varphi(\alpha)$ and $\varphi(\beta)$, by the {\sl Fundamental Theorem of Calculus for the Riemann Integral} we conclude that
$$ \int_{\alpha}^{\beta} f(\varphi(t))\varphi'(t)dt = \int_{\varphi(\alpha)}^{\varphi(\beta)} f(x)dx.$$
The proof is complete. $\qed$

In the following version of the fundamental theorem of calculus for the Lebesgue integral, the function is indeed differentiable at every point belonging to the interval of integration. 

\begin{lemma} \label{TFTCLebesgue} ({\sf The Fundamental Theorem of Calculus for the Lebesgue Integral.}) Let  $f:[a,b]\to \mathbb R$ be differentiable at every point of its domain. Let us suppose that its derivative $f':[a,b]\to \mathbb R$ is Lebesgue integrable. Then, we have the formula
$$ \int_a^b f'(x)dm = f(b) - f(a).$$
\end{lemma}
{\bf Proof} See Rudin [9, pp. 149-150]. $\qed$

\begin{lemma} \label{openness} {\sf (An openness Lemma.)} Consider a continuous map $\varphi: [\alpha,\beta] \to \mathbb R$, differentiable at $p\in (\alpha,\beta)$, with $\varphi'(p)\neq 0$. Then, the following is true.
\begin{itemize}
\item[$\bullet$] There exists $r >0$ such that, for every $|h|<r$ there exists $\kappa$ satisfying
$$ \varphi(p) + h = \varphi(p + \kappa), \ \textrm{with}\ |\kappa| \leq \frac{2}{|\varphi'(p)|}|h|.$$
\end{itemize}
\end{lemma}
{\bf Proof.} Clearly, $\varphi([\alpha,\beta])$ is a bounded, closed, and non-degenerated interval. 

Since $\varphi'(p)\neq 0$, the point $p$ is not a point of minimum, or maximum, of $\varphi$. It thus follows that $\varphi(p)$ is an interior point of $\varphi([\alpha,\beta])$.  It also follows that there exists $\delta_1>0$ such that we have
$$| \varphi(p + \kappa) - \varphi(p) | \geq \frac{|\varphi'(p)|}{2}|\kappa|,\ \textrm{if} \ |\kappa|\leq \delta_1.$$

Given $\epsilon >0$ such that $J= (\varphi(p) - \epsilon, \varphi(p) + \epsilon) \subset \varphi([\alpha,\beta])$, there exists $\delta$, $0<\delta <\delta_1$, satisfying $\varphi\big( (p-\delta,p + \delta) \big) \subset J$. Writing $I = (p-\delta,p + \delta)$, we have
$$\varphi(p) \in \varphi(I) \subset J\subset \varphi([\alpha,\beta]).$$

Since $\varphi'(p)\neq 0$, the point $\varphi(p)$ is not an endpoint of the interval $\varphi(I)$.  Hence, there exists $r>0$ such that 
$$(\varphi(p)- r, \varphi(p) + r) \subset \varphi(I) = \varphi\big((p - \delta, p + \delta)\big).$$
This implies that given $|h|<r$, there exists $|\kappa|<\delta<\delta_1$ such that
$$ \varphi(p) + h = \varphi(p + \kappa), \ \textrm{with} \ |h| \geq \frac{|\varphi'(p)|}{2}|\kappa|.$$
The proof is complete. $\qed$

\begin{lemma} \label{Saks} {\sf (Sak's Theorem on Critical Points.)} Consider a function $\varphi: I \to \mathbb R$, where $I$ is an interval on the real line, and a subset $T\subset I$ such that we have 
$\varphi'=0$ almost everywhere on $T$. Then, $\varphi(T)$ has measure zero.  
\end{lemma}
{\bf Proof.} See Saks [7, p. 227], Serrin and Varberg \cite{Serrin}. The proof is short.

\section{First theorem}

In this section all integrals are Riemann integrals. We employ {\sl Lebesgue's Dominated Convergence Theorem}, but only for a sequence of Riemann integrals. 
We comment more on this in section 10, Remark \ref{R2}. Regarding the substitution maps, we comment on them in Remarks \ref{R3}, \ref{R5} and \ref{R6}. 

\begin{theorem}\label{T1} Let us consider an   integrable $f: I\to \mathbb R$, with $I$ an interval, and a continuous map $\varphi : [\alpha, \beta] \to I$.  Let us suppose that $\varphi$ is differentiable on $(\alpha, \beta)$, with the derivative $\varphi'$ continuous almost everywhere. The following is true.
\begin{itemize}
\item[$\bullet$] If the product $(f\circ\varphi)\varphi'$ is integrable on $[\alpha,\beta]$, then we have the identity
$$\int_{\varphi(\alpha)}^{\varphi(\beta)} f(x)dx =\int_\alpha^\beta f(\varphi(t))\varphi'(t)dt.$$
\end{itemize}
\end{theorem}
{\bf Proof.} Let us split the proof into seven $(7)$ small and easy to follow steps. 
\begin{itemize}
\item[1.] {\sl The setup}.
We may assume that $I$ is the image of $\varphi$. That is, $I=\varphi([\alpha,\beta])$. 

The existence and the value of the Riemann integral of a function do not change if we redefine such function on a finite set of points. Thus, regarding the integral of the product $(f\circ \varphi)\varphi'$ over the interval $[\alpha,\beta]$, if the one-sided derivatives of $\varphi$ at the endpoints $\alpha$ and $\beta$ do not exist, then we may define them at will.

We may assume $f\geq 0$. To show this, it is enough to check that if $f^+$ is the positive part of $f$ then the bounded product function $(f^+
\circ\varphi)\varphi'$ is continuous at every point $t_0$ where $(f\circ \varphi)\varphi'$ and $\varphi'$ are both continuous. The case $\varphi'(t_0)=0$ is trivial since the composite function $(f^+\circ\varphi)$ is bounded. For the case $\varphi'(t_0)\neq 0$, the identity $(f\circ\varphi)  = [(f\circ\varphi)\varphi']/\varphi'$, on a small enough neighborhood of $t_0$, shows that $(f\circ\varphi)$ is continuous at $t_0$ and thus $(f^+\circ\varphi) = \big[|f\circ\varphi| + (f\circ\varphi)\big]/2$ is continuous at $t_0$

Moreover, we assume that the two integrals on the enunciate are finite.
\item[2.] {\sl{Approximating $f$}}. Since the map $\varphi$ is continuous,  its image is a bounded and closed interval. Hence, we may consider a sequence of non-negative functions $f_n : \varphi([\alpha,\beta]) \to [0, + \infty)$ constructed as in section 3. According to Lemma \ref{L1},  each $f_n$ is continuous, satisfies $0\leq f_n \leq f$, and the sequence $(f_n)$ converges pointwise to $f$ at every point of continuity of $f$.

\item[3.]{\sl{ Approximating the product $(f\circ \varphi)\varphi'$}}. Let $\mathcal{N}= \textrm{Discont}(f)$ be the null set of all points of discontinuity of the function $f:\varphi([\alpha,\beta]) \to [0, + \infty)$. Putting $\mathcal{M}= \varphi^{-1}(\mathcal{N})$, we have $\varphi(\mathcal{M})\subset \mathcal{N}$. Thus, the image $\varphi(\mathcal{M})$ is a null set and from Lemma \ref{Serrin} (Serrin and Varberg's theorem) it follows that
$$\left\{\begin{array}{lll}
 \varphi'= 0 \ \textrm{almost everywhere on} \ \mathcal{M},\ \textrm{and}\\
f_n(\varphi(t))\varphi'(t)= f(\varphi(t))\varphi'(t) =0  \ \textrm{almost everywhere on}\ \mathcal{M}.
\end{array}
\right.
$$
Next, let us turn our attention to the outside of the set $\mathcal{M}=\varphi^{-1}(\mathcal{N})$. If a point $t$ is not in $\mathcal{M}$, then  the function $f$ is continuous at $\varphi(t)$ and by applying Lemma \ref{L1} we obtain the convergence $f_n(\varphi(t)) \to f(\varphi(t))$.
Therefore, we have the pointwise convergence
$$f_n(\varphi(t))\varphi'(t) \xrightarrow[\raisebox{0.8ex}{\scriptsize $n\!\to\! +\infty$}]{} f(\varphi(t))\varphi'(t),\ \textrm{almost everywhere on  $[\alpha,\beta]$}.$$
\item[4.] {\sl{ The key integral identity}}. Let $n$ be fixed. By hypothesis, the derivative $\varphi'$ is continuous almost everywhere. Following Lemma \ref{L1}, the function $f_n$ is continuous. These imply that   $(f_n\circ \varphi)\varphi'$ is continuous almost everywhere.  

Yet, again by Lemma \ref{L1}, we know that $0\leq f_n\leq f$. Thus, we have
$$|f_n(\varphi(t))\varphi'(t)| \leq |f(\varphi(t))\varphi'(t)|\ \textrm{on}\  [\alpha,\beta].$$
This inequality shows that the product $(f_n\circ \varphi)\varphi'$ is also bounded on $[\alpha, \beta]$ since, by hypothesis, $(f\circ \varphi)\varphi'$ is Riemann integrable. Summing up, the product $(f_n\circ \varphi)\varphi'$ is Riemann integrable. Even more, since the function $f_n$ is continuous, we know that $f_n$  has a primitive. Hence, we may apply Lemma \ref{Oliveira} (an elementary change of variable for the Riemann integral) and as a result we obtain the integral identity
$$\int_{\varphi(\alpha)}^{\varphi(\beta)} f_n(x)dx = \int_\alpha^\beta f_n(\varphi(t))\varphi'(t)dt.$$
\item[5.] {\sl{ The convergence of the sequence of integrals of $f_n$}}. Lemma \ref{L1} shows that 
$$\int_{\varphi(\alpha)}^{\varphi(\beta)}f_n(x)dx \xrightarrow[\raisebox{0.8ex}{\scriptsize $n\!\to\! +\infty$}]{} \int_{\varphi(\alpha)}^{\varphi(\beta)}f(x)dx.$$
\item[6.] {\sl{The convergence of the sequence of integrals of $(f_{n} \circ \varphi)\varphi'$}}. In step $3$ we showed that  $(f_{n} \circ \varphi)\varphi'$ converges pointwise to $(f\circ \varphi)\varphi'$ almost everywhere. In step 4, that $|(f_n\circ \varphi)\varphi'|$ is Riemann integrable and bounded by the Riemann integrable function $|(f\circ\varphi)\varphi'|$. Therefore, we may employ {\sl Lebesgue's Dominated Convergence Theorem} and then we end up with
$$\int_{\alpha}^{\beta} f_n(\varphi(t))\varphi'(t)dt \xrightarrow[\raisebox{0.8ex}{\scriptsize $n\!\to\! +\infty$}]{} \int_{\alpha}^{\beta} f(\varphi(t))\varphi'(t)dt.$$ 
\item[7.] {\sl{ Conclusion}}. A combination of steps $4, 5$, and $6$ yields the identity
$$\int_{\varphi(\alpha)}^{\varphi(\beta)}f(x)dx = \int_{\alpha}^{\beta}f(\varphi(t))\varphi'(t)dt.$$
\end{itemize}
The proof is complete.$\qed$

\begin{cor} \label{C1} We may suppose that $f$ is zero outside its interval of integration. That is, let $f:I \to \mathbb R$ and $\varphi: [\alpha,\beta] \to I$ be as in Theorem \ref{T1}, including that the product $(f\circ\varphi)\varphi'$ is integrable on $[\alpha,\beta]$. As an additional hypothesis, let us  suppose that the set $\{t\in [\alpha,\beta]: \varphi(t)=\varphi(\alpha) \ \textrm{or}\ \varphi(t)=\varphi(\beta)\}=\varphi^{-1}(\{\varphi(\alpha),\varphi(\beta)\})$ has measure zero. Let us define the function $g:I \to \mathbb R$  as 
$$g=\left\{\begin{array}{ll}
f,\ \textrm{on the closed interval with endpoints}\ \varphi(\alpha)\ \textrm{and} \ \varphi(\beta),\\
0,\ \textrm{elsewhere}.
\end{array}
\right.$$  
Then, we have 
$$\int_{\varphi(\alpha)}^{\varphi(\beta)} g(x)dx = \int_{\varphi(\alpha)}^{\varphi(\beta)} f(x)dx =\int_\alpha^\beta f(\varphi(t))\varphi'(t)=\int_\alpha^\beta g(\varphi(t))\varphi'(t)dt.$$
\end{cor}
{\bf Proof.} We split the proof into four ($4$) steps.
\begin{itemize}
\item[1.] {\sl Notations.} Let $J$ be the closed interval with endpoints $\varphi(\alpha)$ and $\varphi(\beta)$. We consider $\mathcal{X}_J:I\to \mathbb R$, the characteristic function of $J$.
We also consider $\mathcal{X}_{\varphi^{-1}(J)}: [\alpha,\beta]\to \mathbb R$, the characteristic function of $\varphi^{-1}(J)$. 
\item[2.] {\sl The characteristic function $\mathcal{X}_{\varphi^{-1}(J)}: [\alpha,\beta]\to \mathbb R$ is Riemann integrable.} To see this, let us investigate what are the possible points of discontinuity of such bounded function.
We begin by noticing that the substitution  map $\varphi: [\alpha,\beta] \to I$ is continuous, the sets $J\setminus\{\varphi(\alpha),\varphi(\beta)\}$ and $I\setminus J$ are  open on the interval $I$ (i.e., they are relatively open subsets  of $I$), we have the identity $\mathcal{X}_{\varphi^{-1}(J)} \equiv 1$ on the relatively open preimage set $\varphi^{-1}(J\setminus\{\varphi(\alpha),\varphi(\beta)\}) \subset \varphi^{-1}(J)$, and we have the identity $\mathcal{X}_{\varphi^{-1}(J)} \equiv 0$ on the relatively open preimage set $\varphi^{-1}(I \setminus J) = [\alpha,\beta] \setminus \varphi^{-1}(J)$. It thus follows that the set of discontinuities of the characteristic function $\mathcal{X}_{\varphi^{-1}(J)}$ is a subset of $\{t\in [\alpha,\beta]: \varphi(t)= \varphi(\alpha)\ \textrm{or}\ \varphi(t)=\varphi(\beta)\} = \varphi^{-1}(\{\varphi(\alpha),\varphi(\beta\})$, which is supposed to be a set of measure zero. Therefore, we may conclude that the characteristic function $\mathcal{X}_{\varphi^{-1}(J)}$ is Riemann integrable.
\item[3.]{\sl The product $(g\circ\varphi)\varphi'$ is integrable.} We start by noticing that the identity $g= \mathcal{X}_Jf$ and the identity $\mathcal{X}_J(\varphi(t)) = \mathcal{X}_{\varphi^{-1}(J)}(t)$ imply
 that
\begin{eqnarray*}
g(\varphi(t))\varphi'(t) &=& (\mathcal{X}_Jf)(\varphi(t))\varphi'(t)= \mathcal{X}_J(\varphi(t))[f(\varphi(t))\varphi'(t)]\\
                         &=& \mathcal{X}_{\varphi^{-1}(J)}(t)[f(\varphi(t))\varphi'(t)].
\end{eqnarray*}
Hence, from the hypotheses and step 2 it follows that the function $g(\varphi(t))\varphi'(t)$ can be written as a product of two Riemann integrable functions. This shows that $(g\circ\varphi)\varphi'$ is Riemann integrable.
\item[4.]{\sl{Conclusion.}} By step 3, we may apply Theorem \ref{T1} to the function $g: I \to \mathbb R$. We arrive at
$$\int_\alpha^\beta g(\varphi(t))\varphi'(t)dt = \int_{\varphi(\alpha)}^{\varphi(\beta)}g(x)dx = \int_{\varphi(\alpha)}^{\varphi(\beta)}f(x)dx
                                             = \int_\alpha^\beta f(\varphi(t))\varphi'(t)dt.$$
																						The proof is complete.$\qed$
\end{itemize}

\section{Second theorem}

The theorem in this section generalizes Theorem \ref{T1} in the previous section. It's proof is just an adaptation of that  of Theorem \ref{T1} but it relies a little more on {\sl{Lebesgue's Measure and Integration}}. We use the concept of a {\sl{Lebesgue's measurable function}}, an elementary version of {\sl{The Fundamental Theorem of Calculus for the Lebesgue Integral}} (see Lemma \ref{TFTCLebesgue}), and {\sl{ Lebesgue's Dominated Convergence Theorem}}.  The following theorem has no hypothesis on the derivative of the substitution map other than the finiteness of the related integral.

\begin{theorem}\label{T2}
We consider an almost everywhere continuous function $f: I\to \mathbb R$, where $I$ is an interval, a continuous  substitution map $\varphi : [\alpha, \beta] \to I$ that is differentiable on $(\alpha, \beta)$, and the closed interval $J$ with endpoints $\varphi(\alpha)$ and $\varphi(\beta)$. The following is true.
\begin{itemize}
\item[$\bullet$] If $f$ is integrable on $J$, and $(f\circ\varphi)\varphi'$ is integrable on $[\alpha,\beta]$, then we have 
$$\int_{\varphi(\alpha)}^{\varphi(\beta)} f(x)dx =\int_\alpha^\beta f(\varphi(t))\varphi'(t)dt.$$
\end{itemize}
\end{theorem}
{\bf Proof.} We split the proof into seven (7) steps. 
\begin{itemize}
\item[1.] {\sl The setup}. 
As in Theorem \ref{T1}, we may assume that $I$ is the image of $\varphi$ and we may define $\varphi'(\alpha)$ and $\varphi'(\beta)$ at will. We also assume that the two Riemann integrals on the enunciate are finite.

\item[2.] {\sl{Approximating the function $p=f^{+}$, the positive part of $f$}}. We notice that $p$ is almost everywhere continuous on $I=\varphi([\alpha,\beta])$. Let us take the approximating sequence $p_n : \varphi([\alpha,\beta]) \to [0, + \infty)$ given by Lemma \ref{L1}. Each function
$p_n$ is continuous and satisfies $0\leq p_n\leq p $. The sequence $p_n$ converges pointwise to the function $p$ at all points of continuity of $p$.

\item[3.]{\sl{ Approximating the product  $(p\circ \varphi)\varphi'$, a Lebesgue measurable function}}. Echoing the proof of Theorem \ref{T1}, we define $\mathcal{N} = \textrm{Discont}(p)$ as the null set of all points of discontinuity of the function $p:\varphi([\alpha,\beta]) \to [0, + \infty)$. Putting $\mathcal{M}= \varphi^{-1}(\mathcal{N})$,  by employing Lemma \ref{Serrin} (Serrin and Varberg's Theorem) we find that
$$\left\{\begin{array}{lll}
 \varphi'= 0 \ \textrm{almost everywhere on} \ \mathcal{M},\ \textrm{and}\\
p_n(\varphi(t))\varphi'(t)= p(\varphi(t))\varphi'(t) =0  \ \textrm{almost everywhere on}\ \mathcal{M}.
\end{array}
\right.
$$
If $t\notin \mathcal{M}$, then $p$ is continuous at $\varphi(t)$. It does follow the convergence
$$p_n(\varphi(t))\varphi'(t) \xrightarrow[\raisebox{0.8ex}{\scriptsize $n\!\to\! +\infty$}]{} p(\varphi(t))\varphi'(t),\ \textrm{almost everywhere on  $[\alpha,\beta]$}.$$
From Lebesgue measure theory we know that the derivative $\varphi'$ is a Lebesgue measurable function. Hence, since  $p_n\circ \varphi$ is continuous, the product $(p_n\circ\varphi)\varphi'$ is Lebesgue measurable. Thus, the almost everywhere pointwise limit $\lim (p_n\circ\varphi)\varphi' = (p\circ\varphi)\varphi'$ is a Lebesgue measurable function.
\item[4.] {\sl{ The key integral identity.}}
 Let $n$ be fixed. We know that $0\leq p_n\leq p \leq |f|$. Thus, we have the inequality 
$$|p_n(\varphi(t))\varphi'(t)| \leq |f(\varphi(t))\varphi'(t)|\ \textrm{on}\  [\alpha,\beta].$$
This inequality shows that the product $(p_n\circ \varphi)\varphi'$ is bounded on $[\alpha, \beta]$, since $(f\circ \varphi)\varphi'$ is Riemann integrable and thus bounded. It thus follows that the measurable function $(p_n\circ \varphi)\varphi'$ is Lebesgue integrable and that its Lebesgue integral is finite. Even more, since the function $p_n$ is continuous, we know that 
$p_n$  has a primitive $P_n:\varphi([\alpha,\beta]) \to \mathbb R$. Hence, we have the identity $(P_n\circ\varphi)'(t) = p_n(\varphi(t))\varphi'(t)$ at every point $t\in (\alpha,\beta)$. Therefore, we may apply  {\sl{The Fundamental Theorem of Calculus for the Lebesgue Integral}} for a function that is differentiable at every point in its domain and whose derivative is a Lebesgue integrable function (Lemma \ref{TFTCLebesgue}, section 4) to the differentiable composite function $(P_n\circ\varphi): [\alpha + \epsilon, \beta -\epsilon] \to \mathbb R$, for each sufficiently small $\epsilon >0$. We obtain
\begin{eqnarray*}
\int_{\alpha + \epsilon}^{\beta - \epsilon} p_n(\varphi(t))\varphi'(t)dm &=& \int_{\alpha + \epsilon}^{\beta - \epsilon}(P_n\circ\varphi)'(t)dm \\
                                                                         &=& (P_n\circ \varphi)(\beta - \epsilon) - (P_n\circ\varphi)(\alpha + \epsilon).
\end{eqnarray*}
By letting $\epsilon$ goes to zero, we may employ the continuity of the composite function $(P_n\circ\varphi):[\alpha,\beta] \to \mathbb R$ and the fact that the product function $(p_n\circ \varphi)\varphi':[\alpha,\beta] \to \mathbb R$ is a bounded Lebesgue integrable function to conclude that 
$$\int_{\alpha}^{\beta} p_n(\varphi(t))\varphi'(t)dm = P_n(\varphi(\beta)) - P_n(\varphi(\alpha)).$$
Therefore, since $P_n$ is a primitive of $p_n$, by applying  {\sl{The Fundamental Theorem of Calculus for the Riemann Integral}} to the right-hand side of the equation right above we obtain the integral identity
$$ \int_\alpha^\beta p_n(\varphi(t))\varphi'(t)dm = \int_{\varphi(\alpha)}^{\varphi(\beta)} p_n(x)dx .$$
\item[5.] {\sl{ The limit of the sequence of Riemann integrals of $p_n$}}. From Lemma \ref{L1} we obtain the following convergence,
$$\int_{\varphi(\alpha)}^{\varphi(\beta)}p_n(x)dx \xrightarrow[\raisebox{0.8ex}{\scriptsize $n\!\to\! +\infty$}]{} \int_{\varphi(\alpha)}^{\varphi(\beta)}p(x)dx.$$
\item[6.] {\sl{The limit of the sequence of Lebesgue integrals of $(p_n \circ \varphi)\varphi'$}}. In step $3$ we showed that  the product $(p_n\circ \varphi)\varphi'$ converges pointwise to $(p\circ \varphi)\varphi'$, almost everywhere on its domain $[\alpha,\beta]$. In step 4 we showed that the function $|(p_n\circ \varphi)\varphi'|$ is  Lebesgue integrable and bounded by the {\sl Riemann integrable function} 
$|(f\circ\varphi)\varphi'|$. Therefore, we may employ  Lebesgue's Dominated Convergence Theorem and then we end up with
$$\int_{\alpha}^{\beta} p_n(\varphi(t))\varphi'(t)dm \xrightarrow[\raisebox{0.8ex}{\scriptsize $n\!\to\! +\infty$}]{} \int_{\alpha}^{\beta} p(\varphi(t))\varphi'(t)dm.$$ 
\item[7.] {\sl{ Conclusion}}. A combination of steps $4, 5$, and $6$ yields the integral identity
$$\int_{\varphi(\alpha)}^{\varphi(\beta)}p(x)dx = \int_{\alpha}^{\beta}p(\varphi(t))\varphi'(t)dm, \ \textrm{where}\ p = f^+.$$
What we have proven for the function $f$ also holds for the function $-f$.  Moreover,  $f^-$, the negative part of $f$, is equal to the positive part of $-f$. In symbols, $f^-= (-f)^+$. Thus, we have the identity  
$$\int_{\varphi(\alpha)}^{\varphi(\beta)}f^-(x)dx = \int_{\alpha}^{\beta}f^-(\varphi(t))\varphi'(t)dm.$$
Subtracting this identity from the previous one on display, we arrive at
$$\int_{\varphi(\alpha)}^{\varphi(\beta)}f(x)dx = \int_{\alpha}^{\beta}f(\varphi(t))\varphi'(t)dm.$$
Therefore, since $(f\circ \varphi)\varphi'$ is Riemann integrable, we conclude that
$$\int_{\varphi(\alpha)}^{\varphi(\beta)}f(x)dx = \int_{\alpha}^{\beta}f(\varphi(t))\varphi'(t)dx.$$
The proof is complete.$\qed$
\end{itemize}

\begin{cor}\label{C2}
We consider a function $f: I\to \mathbb R$, with $I$ an interval, and a continuous  map $\varphi : [\alpha, \beta] \to I$ that is differentiable on $(\alpha, \beta)$, with $\varphi'$ integrable on each sub-interval $[c,d]\subset (\alpha,\beta)$. Let $f|_J$ be the restriction of $f$ to the closed and bounded interval $J$ with endpoints $\varphi(\alpha)$ and $\varphi(\beta)$. The following is true.
\begin{itemize}
\item[$\bullet$] If $(f\circ\varphi)\varphi'$ is integrable and $f|_J$ is bounded, then $f|_J$ is integrable and  
$$\int_{\varphi(\alpha)}^{\varphi(\beta)} f(x)dx =\int_\alpha^\beta f(\varphi(t))\varphi'(t)dt.$$
\end{itemize}
\end{cor}
{\bf Proof.} We split the proof into three steps.  We may suppose that $I=\varphi([\alpha,\beta]$.

\begin{itemize}
\item[1.] {\sl The setup.} We assume that $(f\circ \varphi)\varphi'$ is integrable on $[\alpha,\beta]$. As in Theorem \ref{T1}, we may define $\varphi'(\alpha)$ and $\varphi'(\beta)$ at will. Then, following Theorem \ref{T2}, it is enough to prove that $f$ is continuous almost everywhere on $I$. 

\item[2.] {\sl About $\varphi,\varphi'$.} Let us fix $K=[c,d]\subset (\alpha,\beta)$. Since $\varphi'$ is integrable on $K$, it follows that $\varphi'$ is almost everywhere continuous and bounded on $K$.

The boundedness of $\varphi'$ on $K$ implies that the restriction $\varphi|_K:K\to \mathbb R$ is a Lipschitz map. We have $|\varphi(t) - \varphi(\tau)|\leq L_K|t- \tau|$, for all $t,\tau \in K$, where $L_K = \sup \{|\varphi'(s)|: s \in K\}$.  Hence, $\varphi|_K$ sends null sets to null sets.

A countable union of null sets is a null set too. Then, it is not hard to see,  $\varphi'$ is continuous a. e.  on $(\alpha,\beta)$, and  $\varphi$ sends null sets to null sets.

\item[3.] {\sl The function $f$ is continuous a.e.} The hypotheses imply that $(f\circ \varphi)\varphi'$ is continuous almost everywhere. Then, by step 2, the set $D$ of all points $t\in (\alpha,\beta)$ such that $\varphi'$ or $(f\circ \varphi)\varphi'$ is discontinuous at $t$  has measure zero. Hence, $\varphi(D)$ is a harmless null subset of $\varphi([\alpha,\beta])$. 

It is enough to analyze the continuity of $f$ on $\varphi(C)$, with $C$ the set of all points $t\in (\alpha,\beta)$ such that $\varphi'$ and $(f\circ \varphi)\varphi'$ are continuous at $t$.

The case $t \in C$, $\varphi'(t)=0$. Then, $t \in Z_{\varphi'}=\{s: \varphi'(s)=0\}$. By Lemma \ref{Saks} (Sak's theorem), the set $\varphi(Z_{\varphi'})$ has measure zero. Thus, 
the image set $\varphi(\{ t\in C: \varphi'(t)=0\})$, a subset of $\varphi(Z_{\varphi'})$, is a null subset of $\varphi([\alpha,\beta])$.

The case $t \in C$ and $\varphi'(t)\neq 0$. By Lemma \ref{openness} (an openness lemma),  there exists $r>0$ such that for every $|h|<r$ there exists $\kappa$ satisfying
$$ \varphi(t) + h = \varphi(t+\kappa), \ \textrm{with} \ |\kappa|\leq \frac{2}{|\varphi'(t)|}|h|.$$
Then,  we conclude the proof of this step by noticing that
$$\lim_{h \to 0} f(\varphi(t) + h)=  \lim_{\kappa \to 0}\frac{\big[(f\circ\varphi)\varphi'\big](t + \kappa)}{\varphi'(t + \kappa)} = \frac{\big[(f\circ\varphi)\varphi'\big](t)}{\varphi'(t)}=f(\varphi(t)) .$$
\end{itemize}
The proof of the corollary is complete. $\qed$

\section{Third theorem}

\begin{theorem}\label{T3} Let us consider an integrable function $f:[a,b]\to \mathbb R$  and a continuous map $\varphi:[\alpha,\beta] \to [a,b]$ that is differentiable on the  open interval $(\alpha, \beta)$. Let us suppose that $\varphi(\alpha) =a$ and $\varphi(\beta)=b$.  The following is true.
\begin{itemize}
\item[$\bullet$] If the product $(f\circ\varphi)\varphi'$ is integrable on $[\alpha,\beta]$, then we have 
$$\int_a^b f(x)dx = \int_{\alpha}^\beta f(\varphi(t))\varphi'(t)dt.$$
\end{itemize}
\end{theorem}
{\bf Proof.} It is immediate from Theorem \ref{T2}. $\qed$

\section{Change of variable for improper integrals}

We say that a function $f: (a,b) \to \mathbb R$ is {\sl improper Riemann integrable}, where  $(a,b)$ is an arbitrary open interval (bounded or unbounded) on the real line, if $f$ is  Riemann integrable on each bounded and closed sub-interval $[u,v]$ of $(a,b)$ and the following limit is finite (the notation for the {\sl improper Riemann integral} of $f$ is as follows),
$$  \lim_{\substack{ u \to a\\ v \to b}}\int_u^v f(x)dx = \int_a^bf(x)dx .$$

\begin{cor} \label{C3}({\sf A Change of Variable Theorem for Improper Riemann Integrals.}) Let us consider a function $f:(a,b)\to \mathbb R$, with $(a,b)$  an open interval on the real line. Let us suppose that either $f:(a,b) \to \mathbb R$ is improper Riemann integrable or $f:[a,b]\to \mathbb R$ is Riemann integrable (supposing that $[a,b]$ is bounded). Let $\varphi:(\alpha,\beta) \to (a,b)$ be a differentiable map, where $(\alpha, \beta)$ is an open interval, that satisfies $\lim_{t\to \alpha}\varphi(t)=a$ and $\lim_{t \to \beta} \varphi(t) =b$.
The following is true.
\begin{itemize}
\item[$\bullet$] If the product $(f\circ\varphi)\varphi':(\alpha, \beta) \to \mathbb R$ is improper Riemann integrable (or Riemann integrable on $[\alpha,\beta]$, if such interval is bounded), then we have
$$\int_a^b f(x)dx = \int_{\alpha}^\beta f(\varphi(t))\varphi'(t)dt.$$
\end{itemize}
\end{cor}
{\bf Proof.} (We notice that the two conditions on possible one-sided limits of $\varphi$ at the values $\alpha$ and $\beta$ [that is, $\varphi(\alpha +) =a$ and $\varphi(\beta -)=b$, if $\alpha$ and $\beta$ are finite] are analogous to the conditions $\varphi(\alpha) =a$ and $\varphi(\beta)=b$ required on Theorem \ref{T3}.)

Let $(\alpha_n)$ and $(\beta_n)$ be sequences on $(\alpha,\beta)$, with $\lim \alpha_n =\alpha$ and $\lim \beta_n = \beta$.  From Theorem \ref{T2} it follows the identity, apropos of two proper Riemann integrals,
$$\int_{\alpha_n}^{\beta_n} f(\varphi(t))\varphi'(t)dt 
                                               =  \int_{\varphi(\alpha_n)}^{\varphi(\beta_n)} f(x)dx.$$ 
To conclude this proof, it is enough to take the limit as $n \to + \infty$ on both sides of the identity right above. The proof is complete. $\qed$

\section{Examples}

\begin{exe} \label{E1} Consider $f(x) = x^3$,  where $x\in (-\infty,+\infty)$, and the map
$$\varphi(t) = t \sin\frac{1}{t}, \ \textrm{if} \ t \in \left(0,\frac{2}{\pi}\right], \ \textrm{and}\ \varphi(0)=0.$$

\end{exe}
\begin{figure}[!ht]
\vspace{-0,4 cm}
\centering
\psfrag{o1a1=o2a2}{\hspace{-0,5 cm}{\fontsize{13}{2}\selectfont$$}}
\psfrag{a}{\hspace{0,0 cm}{\fontsize{15}{2}\selectfont$a$}}
\psfrag{a1}{{\fontsize{15}{2}\selectfont$\cdots$}}
\psfrag{a2}{\hspace{- 0,4 cm}{\fontsize{15}{2}\selectfont$x_{i-1}$}}
\psfrag{dots}{\hspace{- 0,7 cm}{\fontsize{15}{2}\selectfont$$}}
\psfrag{ak1}{\hspace{ - 0,5 cm}{\fontsize{15}{2}\selectfont$\overline{x_i}=\varphi(\overline{t_i})$}}
\psfrag{ak}{\hspace{ -0,0 cm}{\fontsize{15}{2}\selectfont$x_i \ \ \ \ \ \ \ \cdots$}}
\psfrag{b}{\hspace{0,0 cm}{\fontsize{15}{2}\selectfont$b$}}
\resizebox{320pt}{!}{\includegraphics{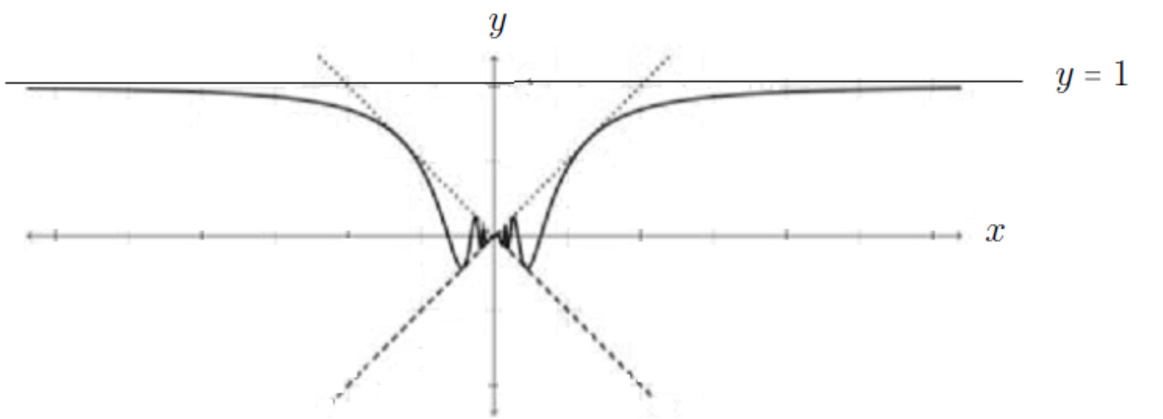}}
\vspace{-0,3 cm}
\caption{A sketch of graph of $\varphi$, supposing $t \in (-\infty,+\infty)$.}
\vspace{-0,0 cm}
\end{figure}

Evidently, the polynomial $f(x)=x^3$ is continuous and integrable on any compact interval, and the  map $\varphi$ is differentiable on $(0,2/\pi)$.  We have
$$\varphi'(t)= \sin\frac{1}{t} - \frac{1}{t}\cos \frac{1}{t},$$ 
which implies that the set of discontinuities of $\varphi'$ is just the single point $x=0$. It also implies that the derivative $\varphi'$ is unbounded and  thus {\sl{not integrable}} on $[0,2/\pi]$.
On other hand, it is not hard to show the integrability of the function
$$(f\circ\varphi)(t)\varphi'(t) = t^3 \left(\sin^3 \frac{1}{t}\right)\left(\sin\frac{1}{t} - \frac{1}{t}\cos\frac{1}{t}\right), \ \textrm{where}\ t\in  [0,2/\pi].$$

Then, we may  apply any one of Theorem \ref{T1}, Theorem \ref{T2} and Theorem \ref{T3}. Whatever our choice, we obtain the integral identity 
$$\int_0^\frac{2}{\pi}x^3\,dx = \int_0^{\frac{2}{\pi}}[\varphi(t)]^3\varphi'(t)\,dt.$$
Thus,
$$\int_0^\frac{2}{\pi}x^3\,dx= \frac{\varphi^4(t)}{4}\Big|_0^{\frac{2}{\pi}}= \frac{4}{\pi^4}.$$
We remark that since the derivative $\varphi'$ is not Riemann integrable, the versions of the Change of Variable Theorem requiring the Riemann integrability of the derivative of the substitution map  do not apply to Example \ref{E1}.

\begin{exe} \label{E2} Consider the function $f(x) = \frac{x}{x^4 +1}$, with $x\in [0,1]$, and  the map $\varphi : [0,1]\to [0,1]$ given by $\varphi(t) = \sqrt{t}$.
\end{exe}

The function $f:[0,1] \to \mathbb R$ is continuous and does integrable. The substitution map $\varphi$ is continuous on $[0,1]$, differentiable on $(0,1]$, and we have
$$\varphi'(t) = \frac{1}{2\sqrt{t}},\ \textrm{for all}\ t \in (0,1].$$
Yet, we notice that the map $\varphi$ satisfies $\varphi(0)=0$ and $\varphi(1)=1$. The product function $(f\circ \varphi)(t)\varphi'(t)$ is given by
$$(f\circ \varphi)(t)\varphi'(t) = \frac{\sqrt{t}}{[(\sqrt{t})^4 + 1]}\frac{1}{2\sqrt{t}} = \frac{1/2}{t^2 +1}, \textrm{for all} t \in (0,1].$$
It does follow that  $(f\circ \varphi)(t)\varphi'(t)$ is integrable on the closed interval $[0,1]$.

Thus, we may apply any of Theorem \ref{T1}, Theorem \ref{T2}, Corolary \ref{C2}, and Theorem \ref{T3} and then write 
$$ \int_0^1 \frac{x}{x^4 +1}dx = \int_0^1 f(\varphi(t))\varphi'(t)dt.$$
By developing the right-hand side of the equation right above, we find that
$$
\int_0^1 \frac{x}{x^4 +1}dx = \frac{1}{2} \int_0^1 \frac{dt}{ t^2 + 1}
														= \frac{\arctan t}{2} \Big|_{t=0}^{t=1} 
														=\frac{\pi}{8}.$$
We remark that since the derivative $\varphi'$  is unbounded (and also continuous)  on the half-open interval $(0,1]$, and thus not Riemann integrable on $[0,1]$, the versions of the Change of Variable Theorem requiring the Riemann integrability of the derivative of the substitution map  do not apply to Example \ref{E2}.

\begin{exe} \label{E3} Consider the function
$f(x) = \frac{1}{x^2 +1},\ \textrm{where}\ x\in (-\infty, + \infty),$
and the substitution map $\varphi: \left(- \frac{\pi}{2}, \frac{\pi}{2}\right) \to (-\infty, + \infty),\ \varphi(\theta) = \tan(\theta)$.
\end{exe}
It is well-known that $f$ is improper Riemann integrable. Clearly, $\varphi$ is differentiable and satisfies the conditions
$\varphi(\theta) \to  \pm \infty$ as  $\theta \to \pm \frac{\pi}{2}$. We also have
$$f (\varphi(\theta))\varphi'(\theta) = \frac{\sec^2 \theta}{ 1 + \tan^2(\theta)} = 1, \ \textrm{for all}\ \theta \in \left(-\frac{\pi}{2}, \frac{\pi}{2}\right).$$
Thus, the product $(f\circ \varphi)\varphi'$ is Riemann integrable on $[-\pi/2, \pi/2]$. Hence, we may apply Corollary \ref{C3} and then we arrive at
$$\int_{-\infty}^{+\infty}\frac{dx}{1 + x^2} = \int_{-\frac{\pi}{2}}^{\frac{\pi}{2}} 1\,d\theta = \pi.$$

\section{Remarks}

\begin{rmk} \label{R1} 
On Lemma \ref{L1}, the choice of the number $\epsilon_n$ was made  to avoid measure theory.  If one assumes that all sub-intervals related to the construction of the function $f_n$ have equal length, then it is possible to write a proof of Lemma \ref{L1}, item 4, that employs Lebesgue's Monotone Convergence Theorem. By considering  the continuous function $g_n=\max(f_1,\ldots, f_n)$, we may suppose that the sequence $(f_n)$ is increasing. It is also possible to proof such item  by using {\sl Osgood's Theorem} and {\sl Arzel\`a's Theorem}, whose proofs employ just a bit of measure theory  (see Thomson [17, pp. 490-493]).
\end{rmk}

\begin{rmk} \label{R2} (Avoiding Lebesgue's Dominated Convergence Theorem.)
To proof Theorem \ref{T1}, step 6, it is possible to avoid Lebesgue's Dominated Convergence Theorem. One can use Lebesgue's Monotone Convergence Theorem instead. To  this effect, one may consider the decomposition of the derivative $\varphi'$ of the substitution map $\varphi$ into its positive and negative parts, $\varphi'= (\varphi')_+ - (\varphi')_-$. Then, by  assuming that $(f_n)$ is an increasing sequence of continuous functions, we deal with the increasing sequences of Riemann integrable functions $(f_n\circ \varphi)(\varphi')_+$ and $(f_n\circ \varphi)(\varphi')_-$. It is also possible to employ {\sl Osgood's Theorem} and {\sl Arzel\`a's Theorem} (see Thomson [17, pp. 490-493]).
\end{rmk}

\begin{rmk} \label{R3} (Substitutions maps.) Theorem \ref{T1} has a very easy corollary. Let us suppose that instead of the hypotheses ``$\varphi'$ almost everywhere continuous'' in Theorem \ref{T1}, we have the hypothesis ``$(f\circ \varphi)$ is almost everywhere continuous and almost everywhere not zero'', keeping all the other hypothesis in the theorem. Under this new set of hypotheses, the derivative $\varphi'$ is almost everywhere continuous. This is true since given a point $t$ such that  the integrable $(f\circ \varphi)\varphi'$ and the composite $(f\circ \varphi)$ are continuous at $t$, with $(f\circ \varphi)(t)\neq 0$, then we may write  
$\varphi' = [(f\circ\varphi)\varphi']/(f\circ\varphi)$ on a small enough neighborhood of $t$.
\end{rmk}

\begin{rmk} \label{R4} (Regarding Kuleshov's remark, \cite{Kuleshov}.)
The enunciate of Theorem \ref{T2} reveals that we may have $f:I\to \mathbb R$ unbounded outside its interval of integration, which is determined by the endpoints $\varphi(\alpha)$ and $\varphi(\beta)$, as long as we keep all other conditions on $f$.
\end{rmk}

\begin{rmk} \label{R5} (Volterra's function.)
As an elementary example of a differentiable map $\varphi:[0,1] \to \mathbb R$ whose derivative $\varphi'$ is bounded but not Riemann integrable, see Goffman \cite{Goffman} and Gordon \cite{Gordon}. Evidently, the set of discontinuities of $\varphi'$ does not have measure zero. The first example of one such map was given in 1881, by V. Volterra \cite{Volterra}. It is important to notice that such example was a striking blow to the intuition of Volterra's contemporaries (see Thomson [17, p. 497]).
\end{rmk}

\begin{rmk} \label{R6} (Substitution maps are not Pompeiu derivatives.) 
There exist strictly increasing bicontinuous maps  $\varphi: [0,1] \to \mathbb [0,1]$ with {\sf bounded} derivatives $\varphi'$ which vanish on a countable dense subset  (see Pompeiu \cite{Pompeiu}, Rooij and Schikhof [11, pp. 80-84], Katznelson and Stromberg \cite{Katznelson}).  It thus follows that these derivatives are positive  and discontinuous on a dense subset.   Let us fix  such a map $\varphi$, and then consider an arbitrary function $f$ such that $f$ and the product function $(f\circ\varphi)\varphi'$ are integrable over $[0,1]$, with $f \neq 0$ almost everywhere. Let us show that we then have $\varphi'=0$ almost everywhere and that this leads to a contradiction. Given a point $t$ such that the almost everywhere continuous  $(f\circ\varphi)\varphi'$ is continuous at $t$, we have the following possibilities. If $f(\varphi(t))=0$, then $\varphi(t)$ is in the null set $f^{-1}(0)$ and, by Lemma \ref{Serrin} (Serrin and Varberg), we have $\varphi'=0$ almost everywhere on $\varphi^{-1}\big(f^{-1}(0)\big)$.  If $f(\varphi(t))\neq 0$, we consider two sub-cases. If $f$ is discontinuous at $\varphi(t)$, then $\varphi(t) \in D_f=\{x: f \ \textrm{ is discontinuous at}\ x\}$, with $D_f$ a null set, and from Lemma \ref{Serrin} it follows that $\varphi'=0$ almost everywhere on $\varphi^{-1}(D_f)$. At last, the case $f(\varphi(t)\neq 0$, with $f$ continuous at $\varphi(t)$. Then, $(f\circ \varphi)$ is continuous and nonzero at $t$, which implies that we may write $\varphi' = \big[(f\circ \varphi)\varphi'\big]/(f\circ \varphi)$, on a small neighborhood of $t$, and conclude that $\varphi'$ is continuous at $t$ and thus, since $\varphi'$ vanishes on a dense set,  we have $\varphi'(t)=0$. All in all, we have $\varphi'=0$ almost everywhere on $[0,1]$. On the other hand, the derivative $\varphi'$ is: defined (finite) for every $t\in [0,1]$, bounded, Lebesgue measurable and Lebesgue integrable. By Lemma \ref{TFTCLebesgue} (the elementary fundamental theorem of calculus for the Lebesgue integral) , the value of such integral is $v= \varphi(1) - \varphi(0) = 1 - 0 = 1$, which shows that the set $\{ t : \varphi'(t)>0\}$ has positive mesure, a contradiction with the fact that we have $\varphi'=0$ almost everywhere,
[An arbitrary real map that is everywhere differentiable on an interval and whose derivative vanishes on a dense subset $D$ of its domain (with the complement of $D$ also dense on its domain), is called a {\sf Pompeiu derivative}.] 

\end{rmk}

\section{Conclusion}

A {\sl rule of thumb} sounds reasonable. 
In practical situations, Theorem \ref{T1} looks good enough. Given a concrete, explicit integrable function $f$, one can pick any differentiable substitution map $\varphi$ that comes up to mind, subject only to the finiteness of the integral of the product $(f\circ\varphi)\varphi'$. This is justifiable since it is highly improbable that one thinks of a $\varphi$ that is a quite hard to build Volterra-type function, or a  $\varphi$  whose derivative is discontinuous on a set of  positive measure. It is also good news that the proof of Theorem 1 is not hard or elaborate, making it suitable for undergraduate courses and for  math-based career professionals. As a final point, the author humbly hopes that this article contributes to further research on the Change of Variable Theorem.

\paragraph{Acknowledgments.}  The author is very grateful to Alberto Torchinsky for his comments and references \cite{Kuleshov} and \cite{Thomson}, and to Alexandre Lymberopoulos, and Nuno Miguel Janu\' ario Alves, for some nice conversations about this important theorem.

\bigskip

\bigskip

\bigskip

\bigskip

\noindent\textit{Departamento de Matem\'atica,
Universidade de S\~ao Paulo\\
Rua do Mat\~ao 1010 - CEP 05508-090\\
S\~ao Paulo, SP - Brasil\\
oliveira@ime.usp.br}

\bigskip


\begin{thebibliography}{99}
\bibitem{Bagby} R. J. Bagby, \textit{The substitution theorem for Riemann integrals}, Real Anal. Exchange, {\bf 27(1)} (2001/2002), 309-314.
\bibitem{Davies} R. O. Davies, \textit{An Elementary Proof of the Theorem on Change of Variable in Riemann Integration}, Math. Gaz., {\bf 45} (1961), 23--25, \\
          \url{http://www.jstor.org/stable/3614765}.
\bibitem{Oliveira} O. R. B. de Oliveira, \textit{Changes of Variable for the Riemann integral on the real line}, Real Anal. Exchange, {\bf 49(2)} (2024), 1-11,\\
         \url{https://doi.org/10.14321/realanalexch.49.2.1689708439}.
\bibitem{Goffman} C. Goffman, \textit{A bounded derivative which is not Riemann integrable}, Amer. Math. Monthly {\bf 84} (1977), 205-206.
\bibitem{Gordon} A. R. Gordon, \textit{A bounded derivative that is not Riemann Integrable}, Math. Mag. {\bf 89(5)} (2016), 364-370,\\
        \url{https://www.jstor.org/stable/10.4169/math.mag.89.5.364}.
\bibitem{Katznelson} Y. Katznelson and K. Stromberg, \textit{Everywhere differentiable, nowhere monotone functions}, Amer. Math. Monthly {\bf 81}  (1974), 349-354.
\bibitem{Kestelman} H. Kestelman,  \textit{Change of variable in Riemann integration}, Math. Gaz., {\bf 45(351)} (1961), 17-23, 
         \url{http://www.jstor.org/stable/3614764}.
\bibitem{Kuleshov} A. Kuleshov, \textit{A Remark on the Change of Variable Theorem for the Riemann integral}, Mathematics, 9, 1899 (2021), \\
         \url{https://doi.org/10.3390/math9161899}.
\bibitem{Pompeiu} D. Pompeiu, \textit{Sur les fonctions d\'eriv\'ees}, Mathematische Annalen (French), {\bf 63 (3)} (1907), 326-332,
        \url{https://doi.org/10.1007%2FBF01449201}.
\bibitem{Pouso} R. L. Pouso, \textit{Riemann integration via primitives for a new proof to the change of variable theorem}, arXiv:1107.1996v1[ [math.CA].
\bibitem{Rooij} A. C. M. van Rooij and W. H. Schikhof, \textit{A Second Course on Real Functions}, Cambridge University Press, 1982.
\bibitem{Rudin} W. Rudin, \textit{Real \& Complex Analysis}. 3rd ed., McGraw-Hill Inc, 1987.
\bibitem{Saks} S. Saks, \textit{Theory of the Integral}, Monografie Matematyczne VII, Warszawa-Lwow, 1937.
\bibitem{Sarkhel} D. N.  Sarkhel and R. V\'yborn\'y,  \textit{A Change of Variables Theorem for the Riemann Integral}, Real Anal. Exch., $\textbf{22}(1)$ (1996-97), 390-395.
\bibitem{Serrin} J. B. Serrin and D. E. Varberg, \textit{A general chain rule for derivatives and the change of variables formula for the Lebesgue integral}, Amer. Math. Monthly {\bf 76} (1969), 514-520.
\bibitem{Tandra} H. Tandra,  \textit{A simple proof of the Change of Variable Theorem for the Riemann Integral}, The Teaching of Mathematics, Vol. {\bf XVIII(1)} (2015), 25-28.
\bibitem{Thomson} B. S. Thomson, \textit{The bounded Convergence Theorem}, Amer. Math. Monthly, {\bf 127(6)} (2020), 483-503, \\
         \url{https://doi.org/10.1080/00029890.2020.1736470}.
\bibitem{Torchinsky1} A. Torchinsky, \textit{The change of variable formulas for Riemann integrals}, Real Anal. Exch., $\textbf{45}(1)$ (2020), 151-172, \\
         \url{https://doi.org/10.14321/realanalexch.45.1.0151}.
\bibitem{Torchinsky2} A. Torchinsky, \textit{A Modern View of the Riemann Integral}, Lecture Notes in Mathematics 2309, Springer, 2022.
\bibitem{Volterra} V. Volterra, \textit{Sui principii del calcolo integrale}, Giorn. Mat. Battaglini, {\bf 19} (1881), 333--372.
\end{thebibliography}
\end{document}